# Aryabhata's Root Extraction Methods

Abhishek Parakh
Louisiana State University
Aug 31st 2006

## 1 Introduction

This article presents an analysis of the root extraction algorithms of Aryabhata given in his book *Āryabhatīya* [1, 2]. Aryabhata, who was born in 476, occupies an important place in the history of mathematics and astronomy. This year's RSA Conference on cryptology honored him, and several studies of his algorithm to solve linear indeterminate equations, which has applications in computer science, have appeared recently [3,4,5].

Aryabhata's book presents his astronomical and mathematical theories. He took the earth to rotate on its axis and he gave planet periods with respect to the sun. The books by Datta and Singh [6] and Srinivasiengar [7] are a good source for a quick history of Indian computing algorithms. For combinatoric and astronomical motivations of Indian mathematics, see the recent papers by Kak [8-11], and for an assessment of Indian computational methods in a wider context, see the book by Joseph [12] and the essay by Pearce [13].

In this paper we analyze Aryabhata's root extraction methods given in the mathematical section of the *Aryabhatiya*, which is called the *Ganitapada*. In his earlier analysis [1], Klintberg showed that Aryabhata's methods were different from that of Greek mathematics. Here we go beyond Klintberg's analysis and discuss the larger question of the theory behind Aryabhata's algorithms. Our analysis shows that Aryabhata was well aware of the place value system of numbers. We also look at the computational complexity of these methods and observe that the methods taught today in schools for root extraction are essentially an extension of Aryabhata's methods.

## 2 The Cube Root extraction method

Aryabhata presents the following concise verse of the description of cube root extraction method:

## GAṆITAPĀDA, VERSE 5

अघनाद् भजेद् द्वितीयात्
त्रिगुणेन घनस्य मूलवर्गेण ।
वर्गस्त्रिपूर्वगुणितः
शोध्यः प्रथमाद् घनश्च घनात् ॥५॥

aghanād bhajed dvitīyāt
triguṇena ghanasya mūlavargeṇa
vargastripūrvaguṇitaḥ
śodhyaḥ prathamād ghanaśca ghanāt[1]

**Translation [2]:** (Having subtracted the greatest possible cube from the last cube place and then having written down the cube root of the number subtracted in the line of the cube root), divide the second non-cube place (standing on the right of the last cube place) by thrice the square of the cube root (already obtained); (then) subtract form the first non-cube place (standing on the right of the second non-cube place) the square of the quotient multiplied by thrice the previous (cube-root); and (then subtract) the cube (of the quotient) from the cube place (standing on the right of the first non-cube place) (and write down the quotient on the right of the previous cube root in the line of the cube root, and treat this as the new cube root. Repeat the process if there is still digits on the right).

**Algorithm:** Represent the given number as a series of indexed digits, i.e., $d_n d_{n-1} \ldots\ldots d_2 d_1 d_0$, where $d_0$ is the digit in units place, $d_1$ is the digit in tens place, $d_2$ is the digit in hundreds place and so on. Let the final root be R and $n$ is the index of the left most digit in the given number. Then the algorithm is as follows:

1. Pick the digit $d_i$ such that $\frac{i}{3}$ = integer,

   such that, $\frac{i+k}{3} \neq$ integer and $\frac{i}{3} > \frac{i-k}{3}$, $k$ =1,2,3, … .

2. $p = \left\lfloor \frac{n-1}{3} \right\rfloor$.

3. Let $k = 100 \times d_{i+2} + 10 \times d_{i+1} + d_i$.
   if $d_{i+2}$ does not exist, let $d_{i+2} = 0$; if $d_{i+1}$ does not exist, let $d_{i+1} = 0$ .

4. Choose $A$ such that $A^3 \le k$ and $k - A^3$ is minimum.

5. $S = k - A^3$.

6. $R = A$.

7. $l = 10 \times S + d_{i-1}$.

8. $S = l \bmod (3 \cdot R^2)$.

9. $m = 10 \times S + d_{i-2}$.

10. $B = \left\lfloor \dfrac{l}{3 \cdot R^2} \right\rfloor$.

11. $S = m - \left(3 \times R \times B^2\right)$.

12. $R = 10 \times R + B$.

13. $n = 10 \times S + d_{i-3}$.

14. $S = n - B^3$.

15. $i = i - 3$.

16. $p = p - 1$.

17. If $p \ne 0$ then go to 7, else quit.

**Example [1]:** Before discussing the theory behind the algorithm, let us first look at its working with the help of an example. Let the given number be $x^3 = 34{,}965{,}783$. We have to estimate the value of $x$.

Step 1: Above each digit, mark out the cubic places (with a dash), and non-cubic places (with a hat), starting from right. The units place is always considered as a cubic place:

$$\bar{3}\ \hat{4}\ \bar{9}\ \bar{6}\ \hat{5}\ \bar{7}\ \bar{8}\ \hat{3}$$

Step 2: Locate the last cubic place (i.e. the second digit from the left) and find the number that holds this position (=34):

$$\begin{array}{cccccccc} & \downarrow & & & & & & \\ \bar{3} & \hat{4} & \bar{9} & \bar{6} & \hat{5} & \bar{7} & \bar{8} & \hat{3} \end{array}$$

Step 3: Now, find the number's nearest lesser (or equal) cube (=27 i.e. $3^3$) and subtract it from 34. Put the cube root of this lesser cube (3) in the root result area. The root result area now contains the most significant digit of the root result:

$$\begin{array}{rcccccccc} & \bar{3} & \hat{4} & \bar{9} & \bar{6} & \hat{5} & \bar{7} & \bar{8} & \hat{3} \\ - & 2 & 7 & & & & & & \\ \hline & & 7 & & & & & & \end{array}$$

Root Result = 3

Step 4: Locate the next place to the right and move its digit (9) down to the right of the result of the previous subtraction to form 79. Since we are now on a second non-cubic place, we divide 79 by three times the square of the current root and evaluate the quotient, i.e. $\left\lfloor \frac{79}{3\times 3^2} \right\rfloor = 2$. Multiply this quotient with three times the square of current assembled root, i.e. $2\times 3\times 3^2 = 54$ and subtract it from 79. This is equivalent to doing a mod operation in the step 8 of our algorithm.

$$\begin{array}{rcccccccc} & \bar{3} & \hat{4} & \bar{9} & \bar{6} & \hat{5} & \bar{7} & \bar{8} & \hat{3} \\ & \ldots\ldots & & & & & & & \\ & & 7 & 9 & & & & & \\ - & & 5 & 4 & & & & & \\ \hline & & 2 & 5 & & & & & \end{array}$$

Root Result = 3

Step 5: Now, bring down the next digit in the number to form 256 and subtract three times the assembled root times the square of the last quotient, i.e. $3\times 3\times 2^2 = 36$. And then place this quotient as the next digit of the assembled root result, i.e. $10\times 3 + 2 = 32$.

$$\begin{array}{rcccccccc} & \bar{3} & \hat{4} & \bar{9} & \bar{6} & \hat{5} & \bar{7} & \bar{8} & \hat{3} \\ & \ldots\ldots & & & & & & & \\ - & & 2 & 5 & 6 & & & & \\ & & & 3 & 6 & & & & \\ \hline & & 2 & 2 & 0 & & & & \end{array}$$

Root Result = 3 2

Step 6: Again, bring down the next digit to the right of the previous subtraction to form 2205 and subtract the cube of the last quotient, i.e. $2^3 = 8$ to get 2197.

$$\bar{3}\ \hat{4}\ \bar{9}\ \bar{6}\ \hat{5}\ \bar{7}\ \bar{8}\ \hat{3}$$

```
 ......
 -2  2 0 5
         8
 ---------
   2 1 9 7
```

Root Result = 3 2

Step 7: Repeat the steps 4, 5 and 6 until the remainder is zero. The resulting calculation is as follows:

$$\bar{3}\ \hat{4}\ \bar{9}\ \bar{6}\ \hat{5}\ \bar{7}\ \bar{8}\ \hat{3}$$

```
 - 2 7
   ---
     7 9
   - 5 4
 ---------
     2 5 6
   -   3 6
 ---------
     2 2 0 5
   -       8
 -----------
     2 1 9 7 7
   - 2 1 5 0 4
 -------------
         4 7 3 8
       - 4 7 0 4
 ---------------
             3 4 3
           - 3 4 3
 -----------------
                 0
```

Root Result = 3 2 7

We see that the procedure is very simple and the intermediate numbers generated are easy to handle computationally. Hence, a large problem is broken down into a series of small steps.

**Theory:** If we look at the binomial expansion $(A+B)^3 = A^3 + 3 \cdot A^2 \cdot B + 3 \cdot A \cdot B^2 + B^3$, the similarity of Aryabhata's method with this expansion is easily noticed. The given number is considered to be $(A+B)^3$ and our aim is to estimate the value of $(A+B)$.

Aryabhata starts out with first subtracting $A^3$ from the given number and the corresponding cube root $A$ becomes the first digit of our final result. The value $A$ is determined by trial and error. If $k$ is the number from which $A$ is to be subtracted, then it is observed in the algorithm that the maximum value of $k$ can be 999. Hence, estimation of $A^3$ to determine $k - A^3$ is very easy and equivalent to one look up table complexity.

In the next step Aryabhata subtracts 3 times the assembled root, which in first iteration is $A$. Then he estimates the value of $B$ and subtracts 3 times the assembled root times the square of the estimate of $B$. In the final step he subtracts the cube of the estimate of $B$. Thus, our first approximation of cube root becomes $10 \cdot A + B$.

After the first iteration of algorithm is done, our first estimate of cube root is obtained. In other words what we have effectively achieved is to subtract the closest cube, in our example $32^3 = 32768$, from the number formed by last 5 digits, that is $34965$.

In each iteration, we improve our knowledge of estimate of cube root by $B$. When we have subtracted exactly $\left(A + B\right)^3$, the final remainder will be zero and the answer we seek will be $\left(A + B\right)$.

## 3 The Square Root extraction method

Aryabhata presents his square root extraction method in the following verse:

GAṆITAPĀDA, VERSE 4

भागं हरेदवर्गान्नित्यं द्विगुणेन वर्गमूलेन ।
वर्गाद्वर्गे शुद्धे लब्धं स्थानान्तरे मूलम् ॥ ४ ॥

bhāgaṁ haredavargānnityaṁ dviguṇena vargamūlena
vargādvarge śuddhe labdhaṁ sthānāntare mūlam[1]

**Translation [2]:** (Having subtracted the greatest possible square from the last odd place and then having written down the square root of the number subtracted in the line of the square root) always divide the even place (standing on the right) by twice the square root. Then, having subtracted the square (of the quotient) from the odd place (standing on the right), set down the quotient at the next place (i.e., on the right of the number already written in the line of the square root). This is the square root. (Repeat the process if there are still digits on the right).

**Algorithm:** Represent the given number as a series of indexed digits, i.e. $d_n d_{n-1} \ldots\ldots d_2 d_1 d_0$, where $d_0$ is the digit in units place, $d_1$ is the digit in tens place, $d_2$ is the digit in hundreds place and so on. Let the final root be R and $n$ is the index of the left most digit in the given number.

1. Pick the digit $d_i$ such that $\frac{i}{2}$ = integer,

   Where, $\frac{i+k}{2} \neq$ integer and $\frac{i}{2} > \frac{i-k}{2}$, $k=1,2,3, \ldots$ .

2. $p = \left\lfloor \frac{n-1}{2} \right\rfloor$.

3. If $d_{i+1}$ exists then $a = 10 \times d_{i+1} + d_i$, else $a = d_i$ .

4. Choose $A$ such that $A^2 \leq a$ and $a - A^2$ is minimum.

5. $S = a - A^2$.

6. $R = A$.

7. $y = 10 \times S + d_{i-1}$.

8. $S = y \bmod (2 \cdot R)$.

9. $B = \left\lfloor \frac{y}{2 \cdot R} \right\rfloor$.

10. $R = 10 \times R + B$.

11. $c = 10 \times S + d_{i-2}$.

12. $S = c - B^2$.

13. $i = i - 2$.

14. $p = p - 1$.

15. If $p \neq 0$ then go to 7, else quit.

**Example [1]:** Let us look at the working of the algorithm. We perform the steps described above using an example: $x^2 = 11943936$ . Our aim is to find the value of $x$ .

Step 1: Mark out the odd places (with a hat) and the even places (with a dash), starting from the right hand side (units place), which is always considered as an odd place:

$$\bar{1}\ \hat{1}\ \bar{9}\ \hat{4}\ \bar{3}\ \hat{9}\ \bar{3}\ \hat{6}$$

Step 2: Locate the last *odd* place (in this case, the second digit from the left) and find the number that holds the position (=11):

$$\begin{array}{cccccccc} & \downarrow & & & & & & \\ \bar{1} & \hat{1} & \bar{9} & \hat{4} & \bar{3} & \hat{9} & \bar{3} & \hat{6} \end{array}$$

Step 3: Find the number's nearest lesser square (in this case $3^2 = 9$) and subtract from the same position. Put its square root (= 3) aside as the most significant digit of the final square root.

$$\begin{array}{cccccccc} \bar{1} & \hat{1} & \bar{9} & \hat{4} & \bar{3} & \hat{9} & \bar{3} & \hat{6} \\ - & 9 & & & & & & \\ \hline & 2 & & & & & & \end{array}$$

Root result: 3

Step 4: Locate next place to the right and move the digit down next to the result of previous subtraction (= 2) to form 29. Since we are in the even place, we divide 29 by twice the current root (= 3). And then put the quotient of $\frac{29}{(2\times 3)}$, i.e. $\left\lfloor \frac{29}{(2\times 3)} \right\rfloor = 4$ to the right of the current root as the second digit. Now, multiply 4 with (2x3) giving 24 and subtract it from 29.

$$\begin{array}{cccccccc} \bar{1} & \hat{1} & \bar{9} & \hat{4} & \bar{3} & \hat{9} & \bar{3} & \hat{6} \\ \ldots & \ldots & & & & & & \\ & 2 & 9 & & & & & \\ - & 2 & 4 & & & & & \\ \hline & & 5 & & & & & \end{array}$$

Root result: 3 4

Step 5: Move the next digit down next to the result of previous subtraction to form 54. Since, we are at an *odd* place; we take the previous quotient and subtract the square of it from 54.

$\bar{1}\ \hat{1}\ \bar{9}\ \hat{4}\ \bar{3}\ \hat{9}\ \bar{3}\ \hat{6}$

………

$$\begin{array}{r} 54 \\ -\ 16 \\ \hline 38 \end{array}$$

Root result: 3 4

Step 6: We repeat the steps 4 and 5 in till we reach the first odd place.

$\bar{1}\ \hat{1}\ \bar{9}\ \hat{4}\ \bar{3}\ \hat{9}\ \bar{3}\ \hat{6}$

$$\begin{array}{l} -\ 9 \\ \hline \ \ 2\ 9 \\ -\ 2\ 4 \\ \hline \ \ \ \ 5\ 4 \\ \ \ -\ 1\ 6 \\ \hline \ \ \ \ 3\ 8\ 3 \\ \ \ -\ 3\ 4\ 0 \\ \hline \ \ \ \ \ \ 4\ 3\ 9 \\ \ \ \ \ -\ \ \ 2\ 5 \\ \hline \ \ \ \ \ \ 4\ 1\ 4\ 3 \\ \ \ \ \ -\ 4\ 1\ 4\ 0 \\ \hline \ \ \ \ \ \ \ \ \ \ \ \ 3\ 6 \\ \ \ \ \ \ \ \ \ \ -\ 3\ 6 \\ \hline \ \ \ \ \ \ \ \ \ \ \ \ \ \ 0 \end{array}$$

Root result: 3 4 5 6

**Theory:** Again it is noticed that binomial expansion $(A+B)^2 = A^2 + 2 \cdot A \cdot B + B^2$ is similar to the steps performed in Aryabhata's method. The given number is considered to be $(A+B)^2$. Our aim is to find $(A+B)$.

Aryabhata starts out with first subtracting $A^2$ from the given number and the corresponding square root $A$ becomes the first digit of our root result. The value of $A^2$ is determined by trial-and-error. However, this step is not computational intensive as there are only 9 numbers with two digit squares. We can set up a look-up table for this operation.

In the next step we subtract 2 times the current assembled root (in this case $A$ ) and estimate $B$ and then subsequently subtract $B^2$ .

Once we have done our first round of subtractions and we are not yet at the last odd position, then we repeat the algorithm again. After the first round of subtractions, we have effectively subtracted $(A+B)^2$ from the given number. This is actually the closest lesser square subtracted from the number formed by the last four most significant digits of the number (in this case, 1194) and our estimated $(A+B)$ = 34 and $(A+B)^2$ =1156. Hence, the present value of the assembled root is taken as $A$ for the next round subtractions.

At every step our knowledge of the root increases by one digit $B$ . When we have reached the final odd place and performed all the subtractions then if the number is a perfect square, we will have a null remainder. This is equivalent to subtracting $(A+B)^2$ .

# 4 A Note on Computational Complexity

The algorithms are observed to be iterative in nature. We first look at the computational complexity of the Cube Root extraction algorithm: step (1) of the algorithm involves one look up table operation and one subtraction; step (2) involves four multiplications, one addition, one division and one subtraction; step (3) again involves four multiplication, one addition, one division and one subtraction; and step (4) again involves the same number of calculations as step (2) and (3). In addition to this, every iteration uses one multiplication and one addition in order to accumulate the root result.

If we represent multiplications by $M$ , additions by $A$ , divisions by $D$ and subtractions by $S$ , then if the number of digits in the given number is $N$ , the number of iterations will be $\frac{N}{3}$ .

The computational complexity of the cube rooting algorithm turns out to be:

$$\frac{N}{3}(\ 12M + 3A + 3D + 4S\ ) + 1 \textit{ look up table operation}$$

In a similar manner the computational complexity of the Square Root extraction algorithm turns out to be:

$$\frac{N}{2}(\ 5M + 3A + 2D + 3S\ ) + 1 \textit{ look up table operation}$$

## 5 Conclusion

We investigated the theory behind Aryabhata's algorithms and found that the methods taught in today's school for root extraction are essentially the same as those presented by Aryabhata centuries ago. Also, the methods are based on the place value system of numbers. The computational complexity of the algorithms was presented.

## 6 Reference

1. B.O. Klintberg, *Were Aryabhata's Square and Cube Root Methods Originally from the Greeks?* M.Sc. Dissertation in Philosophy and Science, London School of Economics and Political Sciences, London, 1998.

2. K.S. Shukla and D.V. Sarma, *Aryabhatiya of Aryabhata.* Indian National Science Academy, 1976.

3. S. Kak, Computational Aspects of the Aryabhata Algorithm, *Indian Journal of History of Science*, 21: 62-71, 1986.

4. T.R.N. Rao and C.-H. Yang, Aryabhata remainder theorem: relevance to crypto-algorithms, *Circuits, Systems, and Signal Processing*, 25: 1-15, 2006.

5. S. Vuppala, The Aryabhata algorithm using least absolute remainders. arXiv: cs.CR/0604012

6. B. Datta and A.N. Singh, *History of Hindu Mathematics, A Source Book*, Parts 1 and 2, (single volume). Asia Publishing House, Bombay, 1962.

7. C.N. Srinivasiengar, *The History of Ancient Indian Mathematics.* The World Press Private, Calcutta, 1967.

8. S. Kak, Aristotle and Gautama on logic and physics. arXiv: physics/0505172

9. S. Kak, The golden mean and the physics of aesthetics. arXiv: physics/0411195

10. S. Kak, Indian Physics: Outline of Early History. arXiv: physics/0310001

11. S. Kak, Birth and Early Development of Indian Astronomy. In *Astronomy Across Cultures: The History of Non-Western Astronomy,* Helaine Selin (ed), Kluwer, pp. 303-340, 2000. arXiv: physics/0101063

12. G. G. Joseph, *The Crest of the Peacock, Non-European Roots of Mathematics*. Princeton University Press, 2000.

13. I.G. Pearce, 2002. Indian mathematics: redressing the balance. http://www-history.mcs.st-andrews.ac.uk/history/Projects/Pearce/index.html